\documentclass[journal]{IEEEtran}
\IEEEoverridecommandlockouts

\def\cdc{0}
\title{
Adding integral action for open-loop exponentially stable
semigroups and application to boundary control of PDE systems
}

\author{A. Terrand-Jeanne, V. Andrieu, V. Dos Santos Martins, C.-Z. Xu\thanks{All authors are with  LAGEP-CNRS, Universit\'e Claude Bernard Lyon1,
Universit\'e de Lyon, Domaine Universitaire de la Doua,
43 bd du 11 Novembre 1918, 69622 Villeurbanne Cedex, France.} }

\usepackage{xcolor}
\usepackage{amsfonts,latexsym,amsmath,amssymb,mathtools}
\usepackage{graphicx,color}
\usepackage{comment}
\usepackage{dsfont}
\usepackage{hyperref}
\newtheorem{theorem}{Theorem}

\newtheorem{proposition}{Proposition}
\newtheorem{corollary}{Corollary}

\newtheorem{assumption}{Assumption}
\newtheorem{remark}{Remark}

\catcode`\@=11
\def\downparenfill{$\m@th\braceld\leaders\vrule\hfill\bracerd$}
\def\overparen#1{\mathop{\vbox{\ialign{##\crcr\crcr
\noalign{\kern0.4ex}
\downparenfill\crcr\noalign{\kern0.4ex\nointerlineskip}
$\hfil\displaystyle{#1}\hfil$\crcr}}}\limits}
\catcode`\@=12

\def\RR{{\mathbb R}}    

\def\BR{\mathcal{B}}    

\def \XR {\mathcal X}
\def \AR {\mathcal A}
\def \BR {\mathcal B}
\def \CR {\mathcal C}
\def \PR {\mathcal P}
\def \MR {\mathcal M}

\DeclareMathOperator{\Id}{I_d}

\def\LR{{\mathfrak L}}  

\def\Un{\mathds{1}}

\begin{document}
\maketitle

\begin{abstract}
The paper deals with output feedback stabilization of exponentially stable systems by an integral controller.
We propose appropriate Lyapunov functionals to prove exponential stability of the closed-loop system.  
An example of parabolic PDE (partial differential equation) systems  and an example of hyperbolic systems are worked out to show how exponentially stabilizing integral controllers are  designed. 
The proof is based on a novel Lyapunov functional construction which employs the forwarding techniques.
\end{abstract}

\section{Introduction}


The use of integral action to achieve output regulation and cancel constant disturbances for infinite dimensional systems has been initiated by S. Pojohlainen in \cite{Pohjolainen_TAC_1982robust}. It has been extended in a series of papers by the same author (see \cite{pohjolainen1985robust} for instance) and some other (see \cite{XuJerbi_IJC_95robust}) always considering bounded control operator and following a spectral approach (see also \cite{paunonen2010internal}).

In the last two decades, Lyapunov approaches have allowed to consider a large class of boundary control problems (see for instance \cite{BastinCoron_2016stability}).
In this work our aim is to follow a Lyapunov approach to solve an output regulation problem. 
The results are separated into two parts. 

In a first part, abstract Cauchy problems are considered. It is shown how a Lyapunov functional can be constructed for a linear system in closed loop with an integral controller when some bounds are assumed on the control operator and for an admissible measurement operator.
This gives an alternative proof to the results of S. Pojohlainen in \cite{Pohjolainen_TAC_1982robust} (and \cite{XuJerbi_IJC_95robust}). 
It allows also to give explicit value to the integral gain that solves the output regulation problem.

In a second part, following the same Lyapunov functional design procedure, we consider a boundary regulation problem for a class of  hyperbolic PDE systems. This result generalizes many others which have been obtained so far in the regulation of PDE hyperbolic systems (see for instance \cite{dos2008boundary,xu2014multivariable,BastinTamasoiu2015stability,TrinhAndrieuXu_TAC_17,TrinhAndrieuXu_Aut_18,TerrandJeanneEtAl_ECC_2018}).

The paper is organized as follows. 
Section \ref{Sec_Abstract} is devoted to the regulation of the measured output for stable abstract Cauchy problems. It  is given a general procedure, for an exponentially stable semigroup in open-loop, to construct a Lyapunov functional for the closed loop system obtained with an integral controller.
Inspired by this procedure, the case of boundary regulation is considered for a general class of  hyperbolic PDE systems in Section \ref{Sec_HypSystems_Theo}. The proof of the  theorem obtained in the context of hyperbolic systems is given in Section \ref{Sec_ProofTheoHyper}.

This paper is an extended version of the paper presented in \cite{TerrandJeanneEtAl_CDC_2018}.
Compare to this preliminary version all proofs are given and moreover more general classes of hyperbolic systems are considered.

\textbf{Notation:} subscripts $t$, $s$, $tt$, $\dots$ denote the first or second derivative w.r.t. the variable $t$ or $s$.
For an integer $n$, $\Id_n$ is the identity matrix in $\RR^{n\times n}$. Given an operator $\AR$ over a Hilbert space, $\AR^*$ denotes the adjoint operator. $\mathcal D_n$ is the set of diagonal matrix in $\RR^{n\times n}$.

\section{General Abstract Cauchy problems}
\label{Sec_Abstract}
\subsection{Problem description}

Let $\XR$ be a Hilbert space with scalar product $\langle,\rangle_\XR$ and $\AR:D(\AR)\subset\XR\rightarrow \XR$  be the infinitesimal generator of a $C_0$-semigroup denoted $t\mapsto e^{\AR t}$.
Let $\BR$ and $\CR$ be  linear operators, $\BR$ from $\RR^m$ to $\XR$ and $\CR$ from $D(\CR)\subseteq\XR$ to $\RR^m$.

In this section, we consider the controlled Cauchy problem with output $\Sigma(\AR,\BR,\CR)$ in Kalman form, as follows
\begin{equation}\label{eq_AbsCauchyProblem}
\varphi_t = \AR \varphi + \BR u + w\ ,\ y=\CR \varphi ,
\end{equation}
where $w \in \XR$ is an unknown constant vector and $u:\RR_+\mapsto\RR^m$ is the controlled input.
We consider the following exponential stability property for the operator $\AR$.
\begin{assumption}[Exponential Stability]\label{ass_EXPStab}
The operator $\AR$ generates a $C_0$-semigroup  which is exponentially stable.
In other words, there exist $\nu$ and $k$ both positive constants such that,  $\forall \varphi_0\in \XR$ and $t \in \RR_+$
\begin{equation}\label{eq_ExpStab}
\|e^{\AR t}\varphi_0\|_\XR \leqslant   k\exp(-\nu t)\|\varphi_0\|_\XR.
\end{equation}
\end{assumption}

\par\vspace{0.5em}
We are interested in the regulation problem. More precisely we are concerned with the problem of regulation of the output $y$ via the integral control
\begin{equation}\label{eq_ControlInte}
u =  k_i K_i z \ ,\ z_t = y -y_{ref},
\end{equation}
where $y_{ref} \in \RR^m$ is a prescribed reference,
 $z \in \RR^m$,  $K_i \in \RR^{m\times m}$ is a full rank matrix  and $k_i$  a positive real number.

The control law being dynamical, the state space has been extended. 
Considering the system $\Sigma(\AR,\BR,\CR)$ in closed loop, with integral control law given in (\ref{eq_ControlInte}) the state space is now $\XR_e = \XR\times\RR^m$ which is a Hilbert space with inner product 
$$
\langle \varphi_{ea},\varphi_{eb}\rangle_{\XR_e} = \langle \varphi_{a},\varphi_{b}\rangle_\XR + z_a^\top z_b \ ,
$$ where  $\varphi_{ea} = \begin{bmatrix} \varphi_a \\ z_a\end{bmatrix}$ and $\varphi_{eb} = \begin{bmatrix} \varphi_b\\z_b\end{bmatrix}$. The associated  norm is denoted $\|\cdot\|_{\XR_e}$.
Let $\AR_e: (D(\AR)\cap D(\CR) )\times\RR^m  \rightarrow \XR_e$ be the extended operator defined as
\begin{equation}\label{def_Ae}
\AR_e  = 
\begin{bmatrix}
\AR & \BR K_i k_i \\
\CR & 0
\end{bmatrix}\ .
\end{equation}

The regulation problem to solve can be rephrased as the following.

\noindent\textbf{Regulation problem:} We wish to find a positive real number $k_i$ and a full rank matrix $K_i$ such that $ \forall (w,y_{ref}) \in \XR\times\RR^m$:
\begin{enumerate}
\item The system \eqref{eq_AbsCauchyProblem}-\eqref{eq_ControlInte} is well-posed. In other words, for all $ \varphi_{e0}=( \varphi_0,z_0) \in \XR_e$ there exists a unique (weak) solution denoted $ \varphi_e(t)=\begin{bmatrix}\varphi(t)\\z(t)\end{bmatrix} \in C^0(\RR_+,\XR_e)$ defined $\forall t\geqslant   0$ and initial condition $ \varphi_e(0)= \varphi_{e0}$.
\item There exists an equilibrium point denoted $ \varphi_{e\infty}=\begin{bmatrix}  \varphi_{\infty} \\ z_\infty\end{bmatrix} \in \XR_e$, depending on $w$ and $y_{ref}$, which is exponentially stable for the system \eqref{eq_AbsCauchyProblem}-\eqref{eq_ControlInte}. In other words, there exist positive real numbers $\nu_e$ and $k_e$ such that for all $t\geqslant   0$
$$
\|\varphi_e(t)-\varphi_{e\infty}\|_{\XR_e}\leqslant   k_e\exp^{-\nu_e t}\|\varphi_{e0}-\varphi_{e\infty}\|_{\XR_e}.
$$
\item The output $y$ is regulated toward the reference $y_{ref}$. More precisely,  
\begin{equation}\label{eq_RegulDef}
\forall \varphi_{e0},\,\, \lim_{t\rightarrow +\infty} |\CR \varphi(t)-y_{ref}|=0 .
\end{equation}
\end{enumerate}

\par\vspace{0.5em}
We know with the work of S. Pohjolainen in \cite{pohjolainen1985robust}
that when $\AR$ generates an exponentially stable analytic semi-group, when $\BR$ is bounded and when $\CR$ is $\AR$-bounded, with a rank condition, the regulation may be achieved.
This result has been extended to more general exponentially stable semi-groups in \cite{XuJerbi_IJC_95robust}.
\begin{theorem}[\cite{XuJerbi_IJC_95robust}]
\label{Theo_XuJerbi}
Assume that $\XR$ is separable and that $\AR$ satisfies Assumption \ref{ass_EXPStab}.
Assume moreover that~:
\begin{enumerate}
\item the operator $\BR$ is bounded;
\item the operator $\CR$ is $\AR$-admissible (see \cite{tucsnak2009observation}), i.e.
\begin{itemize} 
\item it is $\AR$-bounded~:
\begin{equation}\label{eq_ABounded}
|\CR \varphi| \leqslant   c(\|\varphi\|_\XR + \|\AR \varphi\|_\XR)\ ,\ \forall \ \varphi\in D(\AR),
\end{equation}
for some positive real number $c$;
\item  there exist $T>0$ and $c_T>0$ such that 
$$
\int_0^T |\CR e^{\AR t}\varphi|^2dt \leqslant    c_T^2 \|\varphi\|_\XR^2\ ,\ \forall \ \varphi\in D(\AR);
$$
\end{itemize}
\item the rank condition holds. In other words operators $\AR$, $\BR$ and  $\CR$ satisfy
\begin{equation}\label{eq_RankCondition}
\texttt{rank}\{\CR \AR^{-1}\BR\} = m;
\end{equation}
\end{enumerate}
then there exists  a positive real number $k_i^*$ and a $m\times m$ matrix $K_i$, such that for all $0<k_i<   k_i^* $ the operator $\AR_e$ given in \eqref{def_Ae}
is the generator of an exponentially stable $C_0$-semigroup in the extended state space $\XR_e$. More precisely, the system (\ref{eq_AbsCauchyProblem}) in closed loop with (\ref{eq_ControlInte}) is well-posed and the equilibrium is exponentially stable. Moreover,  for all $w$ and $y_{ref}$, equation \eqref{eq_RegulDef} holds (i.e the regulation is achieved).
\end{theorem}

\par\vspace{0.5em}
On another hand, if one wants to address nonlinear abstract Cauchy problems or unbounded operators, we may need to follow a Lyapunov approach.
For instance in the context of boundary control, a Lyapunov functional approach has allowed to tackle feedback stabilization of a large class of PDEs (see for instance \cite{BastinCoron_2016stability} or \cite{coron2007control}).

It is well known (see for instance \cite[Theorem 8.1.3]{JacobZwart_Book_12}) that exponential stability of 
the operator $\AR$ is equivalent to existence 	
of a bounded positive and self adjoint operator $\PR$ in $\LR(X)$ such that 
\begin{equation}\label{LyapEq}
\langle \AR \varphi,\PR \varphi\rangle_\XR  + \langle \PR \varphi, \AR \varphi\rangle_\XR \leqslant   -\mu\|\varphi\|_\XR \ , \ \forall \ \varphi\in D(\AR),
\end{equation}
where $\mu$ is a positive real number.
We assume that this Lyapunov operator $\PR$ is given.
The first question, we intend to solve is the following:
\textit{Knowing the Lyapunov operator $\PR$, is it possible to construct a Lyapunov operator $\PR_e$ associated to the extended operator $\AR_e$?}

To answer this question, we first give a construction based on a well-known technique  in the nonlinear finite dimensional control community named \textit{the forwarding} (see for instance \cite{MazencPraly_TAC_96}, \cite{SepulchreJankovicKokotovic_Aut_97} or more recently \cite{BenachourEtAl_TAC_13}, or \cite{astolfi2017integral}).

\subsection{A Lyapunov approach for regulation}
\label{Sec_LyapAbstract}
Inspired by the forwarding techniques, the following result can be obtained.
\begin{theorem}[Forwarding Lyapunov functional]
\label{Theo_Abstract} Assume that all assumptions of Theorem \ref{Theo_XuJerbi} are satisfied
and let $\PR$ in $\LR(X)$ be a positive self adjoint operator  such that 
(\ref{LyapEq}) holds.
Then there exist a bounded operator $\MR:\XR\rightarrow \RR^m$ and positive real numbers $p$ and $k_i^*$, such that for all $0<k_i<   k_i^*$, there exists $\mu_e>0$ such that the operator
\begin{equation}
\PR_e = \begin{bmatrix}
 \PR + p\MR^*\MR & -p \MR^* \\-p\MR &p \Id
\end{bmatrix}
\end{equation}
is positive and satisfies $\forall \ \varphi_e = (\varphi,z)\in D(\AR)\times \RR^m$
\begin{multline}
\hspace*{-0.35cm}\langle \AR_e\varphi_e,\PR_e\varphi_e\rangle_{\XR_e}  + \langle \PR_e\varphi_e,\AR_e\varphi_e\rangle_{\XR_e} 
\leqslant   -\mu_e( \|\varphi\|_\XR^2 + |z|^2).
\end{multline}
\end{theorem}

\par\vspace{0.5em}
\begin{proof}  
The operator $\AR$ satisfying Assumption \ref{ass_EXPStab}, $0$ is in its resolvent set and consequently $\AR^{-1}:\XR\mapsto D(\AR)$ is well defined and bounded.
Let $\MR:\XR\rightarrow\RR^m$ be defined by $\MR=\CR \AR^{-1}$ which is well defined due to the fact that $D(\AR)\subseteq D(\CR)$ since $\CR$ is $\AR$-bounded.
Moreover, with (\ref{eq_ABounded})  $\forall\varphi \in \XR$ 
$$
|\MR \varphi| = |\CR \AR^{-1}\varphi|\leqslant   c\left( \|\AR^{-1}\varphi\|_\XR + \|\varphi\|_\XR\right) \leqslant   \tilde c \|\varphi\|_\XR,
$$
where $\tilde c$ is a positive real number.
Hence, $\MR$ is a bounded linear operator.
Moreover, $\MR$ satisfies the following equation
\begin{equation}\label{eq_M}
\MR \AR \varphi = \CR \varphi\ ,\ \forall \varphi\in D(\AR)\ .
\end{equation}
Let $K_i = (\CR \AR^{-1}\BR)^{-1}$ which exists due to the third assumption of Theorem 1.
Note that,
\begin{equation}
\langle \varphi_e,\PR_e\varphi_e\rangle_{\XR_e} = \langle \varphi,\PR \varphi\rangle_{\XR}  +  p (z-\MR \varphi)^\top (z-\MR \varphi),
\end{equation}
hence $\PR_e$ is positive.
This candidate Lyapunov functional is similar to the one given in \cite[Equation (34)]{BenachourEtAl_TAC_13}. It is selected following a forwarding approach.

Moreover, we have
\begin{multline*}
\langle \AR_e \varphi_e,\PR_e \varphi_e\rangle_{\XR_e}  + \langle \PR_e \varphi_e,\AR_e \varphi_e\rangle_{\XR_e} =
\\\langle \AR \varphi, \PR \varphi\rangle_{\XR} +\langle \PR \varphi, \AR \varphi\rangle_{\XR}  
\\+ 2p (z-\MR \varphi)^\top (\CR \varphi-\MR \AR \varphi) + k_i\langle \varphi,\PR \BR K_i z \rangle_\XR \\
+  k_i\langle \PR \BR K_i z,\varphi \rangle_\XR -2p(z-\MR \varphi)^\top \MR \BR K_ik_i z.
\end{multline*}
Employing equation (\ref{eq_M}) and $\MR \BR K_i=\Id_m$, the former inequality becomes
\begin{multline}\label{eq_Truc}
\langle \AR_e \varphi_e,\PR_e \varphi_e\rangle_{\XR_e}  + \langle \PR_e \varphi_e,\AR_e \varphi_e\rangle_{\XR_e} =\\
\langle \AR \varphi, \PR \varphi\rangle_{\XR} +\langle \PR \varphi, \AR \varphi\rangle_{\XR}  
+ k_i\langle \varphi,\PR \BR K_i z \rangle_\XR 
\\+ k_i\langle \PR \BR K_i z,\varphi \rangle_\XR -2p(z-\MR \varphi)^\top k_i z.
\end{multline}
Let $\|\PR \BR K_i\|_\XR^2=\alpha$ which is well defined due to the boundedness assumption on $\BR$.
Given $ a, b$ positive constants, the following inequalities hold
\begin{align}\label{eq_KeyStep}
\langle \varphi,\PR \BR K_i z \rangle_\XR &\leqslant   \frac{1}{2a}\|\varphi\|_\XR^2 + \frac{a \alpha}{2}|z|^2,\\
z^\top \MR \varphi &\leqslant   \frac{1}{2b}\|\varphi\|^2 + \frac{b\|\MR\|^2}{2}|z|^2,
\end{align}
it yields given (\ref{LyapEq}) that
\begin{multline}
\langle \AR_e \varphi_e,\PR_e \varphi_e\rangle_{\XR_e}  + \langle \PR_e \varphi_e,\AR_e\varphi_e\rangle_{\XR_e} 
\\
\leqslant   \left[-\mu+\frac{k_i}{a}+\frac{pk_i}{b}\right]  \|\varphi\|_\XR^2
\\+k_i\left( p(-2 + b\|\MR\|^2) + a\alpha  \right)|z|^2 .
\end{multline}
We pick $b$ sufficiently small such that
\begin{equation}\label{eq_Step1}
-2 + b\|\MR\|^2 <0 .
\end{equation}
In a second step, we select $a$  sufficiently small and $p$ sufficiently large such that
\begin{equation}\label{eq_Step2}
p(-2 + b\|\MR\|^2) + a \alpha <0 .
\end{equation}
Finally, picking $k_i^*$ sufficiently small such that
\begin{equation}\label{eq_Step3}
 -\mu+\frac{k_i^*}{a}+\frac{pk_i^*}{b} <0
\end{equation}
 the result is obtained with \\
$ \hspace*{0.5cm}
\mu_e = \min\left\{\mu-\frac{k_i}{a}-\frac{pk_i}{b},p(2 - b\|\MR\|^2) - a \alpha\right\}.
$
\end{proof}
\subsection{Discussion on the result}
A direct interest of the Lyapunov approach given in Theorem \ref{Theo_Abstract}, is that it allows to give an explicit value for $k_i^*$ which appears in Theorem \ref{Theo_XuJerbi}. 
\if\cdc1
\else
We may compute the largest value of $k_i^*$ following this route.
First of all, from (\ref{eq_Step3})
\begin{eqnarray}
k_i^*&=&  \sup_{a,b,p, \text{such that }(\ref{eq_Step1})-(\ref{eq_Step2})}\left\{\frac{\mu p}{\frac{p}{a}+\frac{p^2}{b}}\right\},\\
&=&  \mu\sup_{a,b,p, \text{such that }(\ref{eq_Step1})-(\ref{eq_Step2})}
\left\{\frac{a b}{{p}{a}+{b}}\right\}
\end{eqnarray}
On another hand, taking the the value of $a$ and $b$ given by (\ref{eq_Step1}) and (\ref{eq_Step2}), one can  rewritten them with $ 0<\beta<1$ and $0<\theta<1$ as
\begin{eqnarray}\label{eq_Step4}
b=\frac{ 2}{\|\MR\|^2}\beta,\,\,
 a= 2(1-\beta) \theta \frac{p}{\alpha }.
\end{eqnarray}
Then \begin{eqnarray}
\frac{a b}{{p}{a}+{b}}=\frac{2(1-\beta)\beta  \theta {p} }{{(1-\beta) \theta {\|\MR\|^2}{p^2}}+{{ \alpha\beta}}}
\end{eqnarray}
which right expression is a function of $p$ taking its maximum value when $$p=\sqrt{\frac{\alpha\beta}{(1-\beta) \theta {\|\MR\|^2}}}$$ then
\begin{eqnarray}
k_i^*&=&  \sup_{a,b,p, \text{such that }(\ref{eq_Step1})-(\ref{eq_Step2})}\left\{\frac{\mu p}{\frac{p}{a}+\frac{p^2}{b}}\right\},\\
&=&  \mu\sup_{0<\theta<1,\,\, 0<\beta<1}
\left\{\frac{\sqrt{\beta(1-\beta)\theta}}{\sqrt{\alpha}\|\MR\|}\right\}
\end{eqnarray}
It is reached for $\beta=\frac{1}{2}$ and $\theta =1$ and this yields
$$
k_i^* = \frac{\mu }{2\|\MR\| \sqrt{\alpha}}= \frac{\mu}{2\|\CR \AR^{-1}\| \|\PR \BR (\CR \AR^{-1}\BR)^{-1}\|_\XR} .
$$

Of course, this optimal value depends on the considered Lyapunov operator $\PR$ solution of (\ref{LyapEq}). 
Note that a possible solution to this equation with $\mu=1$ is given for all $(\varphi_1,\varphi_2)$ in $\XR^2$ by (see \cite{JacobZwart_Book_12}) 
$$
\langle \varphi_1, \PR \varphi_2\rangle_\XR = \lim_{t\rightarrow +\infty} \int_0^t \langle e^{\AR s}\varphi_1, e^{\AR s}\varphi_2\rangle_\XR ds .
$$
Due to (\ref{eq_ExpStab}) it is well defined and positive. 
Note also that we have 
$$
\|\PR \|_\XR \leqslant   \frac{k^2}{2\nu}.
$$
This implies, the following corollary.
\fi
\begin{corollary}[Explicit integral gain]\label{Corollary_ExplicitGain}
Given a system $\Sigma(\AR,\BR,\CR)$ satisfying the assumptions of the Theorem \ref{Theo_XuJerbi}, points 1), 2), and 3) of Theorem \ref{Theo_XuJerbi} hold with $K_i = (\CR\AR^{-1}\BR)^{-1}$ and 
\begin{equation}\label{eq_Explicitgain}
k_i^*=\frac{\nu}{\|\CR \AR^{-1}\| k^2 \|\BR (\CR \AR^{-1}\BR)^{-1}\|}.
\end{equation}
\end{corollary}
\vspace*{0.2cm}
An interesting question would now to know in which aspect this value may be optimal.
\vspace*{0.2cm}

\subsection{Illustration on a parabolic systems}
Consider the problem of heating a bar of length $L=10$ with both endpoints at temperature zero. We control the heat flow in and out  around the points $s=2,\; 5,\;$ and $7$
and measure the temperature at points $3,\;6,\;$ and $8$. The problem is to find an integral controller such that the measurements at $s=3,\; 6,\;$ and $8$ are regulated to (for instance) $1$, $3$, and $2$, respectively.  Thus the control system is governed by the following PDE
\begin{multline}
\phi_t(s,t) =
\phi_{ss}(s,t)+ \Un_{[\frac{3}{2}, \frac{5}{2}]}(s) u_1(t) + \Un_{[\frac{9}{2}, \frac{11}{2}]}(s) u_2(t)
\\
+
\Un_{[\frac{13}{2}, \frac{15}{2}]}(s) u_3(t),\;\; (s,t) \in (0, \;10) \times (0, \infty)\end{multline}
where $\phi:[0,+\infty)\times[0,10]\rightarrow \RR$ with boundary conditions
$$\phi(0,t)=\phi(10,t)=0$$
\begin{equation} \label{parabolic1} \phi(s,0)= \phi_0(s),\end{equation}
where $\Un_{[a, b]}:[0,10]\rightarrow \RR$ denotes the characteristic function on the interval $[a, b]$, i.e.,
$$\Un_{[a, b]}(s)=\left\{\begin{array}{ll}1 & \forall\; s\in [a,b],\\
0 & \forall\; s \not \in [a,b].\end{array}\right. $$
The output and the reference are given as
$$y(t)=\begin{bmatrix}\phi(t,3)\\
\phi(t,6)\\ \phi(t,8)\end{bmatrix}\ ,\ y_{ref}=\begin{bmatrix} 1\\ 3\\ 2\end{bmatrix}.
$$

Let the state space be the Hilbert space $\mathcal{X} = L^2((0,10),\RR)$ with usual inner product, and let
the input space and the output space be equal to $\mathbb R^3$. Clearly, from (\ref{parabolic1}), we get
the semigroup generator $\mathcal{A}: D (\mathcal{A}) \rightarrow \mathcal{X}$, the input operator
$\mathcal{B}: \mathbb R^3 \rightarrow \mathcal{X}$ and the output operator
$\mathcal{C}: D (\mathcal{A}) \rightarrow \mathbb R^3$ as follows:
$$ D (\mathcal{A}) =\{\varphi\in H^2(0,10)\;\vert \; \varphi(0)=\varphi(10)=0\},$$
and
$$ \mathcal{A}\varphi=\varphi_{ss}\;\;\forall\; \varphi\in D (\mathcal{A}),$$
$$ \mathcal{B} u= \Un_{[\frac{3}{2}, \frac{5}{2}]} u_1 +\Un_{[\frac{9}{2}, \frac{11}{2}]} u_2+
\Un_{[\frac{13}{2}, \frac{15}{2}]} u_3,$$
and
$$
C \varphi= \begin{bmatrix}\varphi(3)\\ \varphi(6)\\\varphi(8)\end{bmatrix}.
$$
Moreover, note that with Sobolev embedding, an integration by part and by completing the square, 
we have for all 
$\varphi$ in $D(\AR)$
\begin{align*}
\sup_{s\in (0,10)}|\varphi(s)|
&\leq c\int_0^{10} \varphi(s)^2ds+c\int_0^{10} \varphi_s(s)^2ds\\
&\leq c \|\varphi\|_{\XR}+c\int_0^{10} |\varphi(s)\varphi_{ss}(s)|ds\\
&\leq \frac{3}{2}c \|\varphi\|_{\XR} +  \frac{1}{2}c\|\varphi_{ss}\|_{\XR}.
\end{align*}
Hence $\CR$ is $\AR$-bounded. 

Moreover, by direct computation we find that
$$ \mathcal{C} \mathcal{A}^{-1} \mathcal{B}= \frac{-1}{10}
 \begin{bmatrix}14 & 15 &  9\\ 8 & 20 & 18 \\ 4 & 10 & 14\end{bmatrix}.$$
 It is easy to see that the above  matrix is regular.
 Consequently all Assumptions of Theorem \ref{Theo_XuJerbi} hold. With Corollary \ref{Corollary_ExplicitGain}, it is possible to compute explicitly the integral controller gain.
 By direct computation we have for all $\varphi$ in $\XR$
$$
\mathcal{C} \mathcal{A}^{-1} \varphi =
\begin{bmatrix}\frac{3}{10}\int_0^{10} (s-10) \varphi(s) d s  + \! \int_0^3 (3-s) \varphi(s) d s  \\[0.2cm]
\frac{3}{5}\int_0^{10} (s-10) \varphi(s) d s +  \! \int_0^6 (6-s) \varphi(s) d s \\[0.2cm]
\frac{4}{5}\int_0^{10} (s-10) \varphi(s) d s  + \!  \int_0^8 (8-s) \varphi(s) d s \end{bmatrix},
$$
which gives $\| \mathcal{C} \mathcal{A}^{-1}\| \leq 6.2466$.
We have
 $$K_i=  \begin{bmatrix} -1.250 &    1.500 &  -1.125\\
    0.500  & -2.000   &  2.250\\
         0 &    1.000  & -2.000 \end{bmatrix}.$$
 For the open-loop system, consider the Lyapunov operator $\mathcal{P}=\Id$. Then the growth rate
 may be taken as $\mu=\frac{\pi^2}{50}$. 
 It is easy to see that 
 $\| K_i \| =4.2433$, and $\|\mathcal{B}\|\leq \sqrt{3}$. Putting together the numerical values into the
 formula \eqref{eq_Explicitgain} allows to estimate the tuning parameter
 $$k_i^* = \frac{\textstyle \omega }{\textstyle 2
 \|\mathcal{B} K_i\| \;  \| \mathcal{C} \mathcal{A}^{-1}\|} \approx   2.1498*10^{-3}.$$
 With Corollary \ref{Corollary_ExplicitGain}, the integral controller \eqref{eq_ControlInte} with $ 0< k_i< 2.1498*10^{-3}$  stabilizes
 exponentially the equilibrium along  solutions of the closed-loop system and drives asymptotically the measured temperatures to the reference values for any initial condition.

\section{Case of boundary regulation for hyperbollic PDEs}
\label{Sec_HypSystems_Theo}

In the following section we adapt this framework to hyperbolic PDE systems with boundary control.

\subsection{System description}
To illustrate the former abstract theory, we consider the case of hyperbolic partial differential equations as studied in \cite{CoronBastinAndreaNovel_2008dissipative}.
More precisely, the system is given by a one dimensional $n\times n$ hyperbolic system
\begin{multline}
\label{eq_Hyp}
\phi_t(s,t) + \Lambda_0(s) \phi_s(s,t)  +  \Lambda_1(s)\phi(s,t)=0 \\ s\in(0,1),\ t\in[0,+\infty),
\end{multline}
where $\phi:[0,+\infty)\times [0,1] \rightarrow \RR^n$
\begin{align*}
\Lambda_0(s) &= \texttt{diag}\{ \lambda_1(s),\dots,\lambda_n(s)\}\\
\lambda_i(s) &>0\ \forall i\in\{1,\dots,\ell\}\\
\lambda_i(s) &<0\ \forall i\in\{\ell + 1,\dots,n\}
,
\end{align*}
where the maps $\Lambda_0$ is in $C^1([0,1];{\mathcal D}_n)$ and $\Lambda_1$ is in $C^1([0,1];\RR^{n\times n})$  with the initial condition $\phi(0,s) = \phi_0(s)$ for $s$ in $[0,1]$
where $\phi_0:[0,1] \rightarrow \RR^n$
and with the boundary conditions
\begin{eqnarray}\label{Hyp_BC}
&&\hspace*{-0.75cm}\begin{bmatrix}
\phi_+(t,0)\\
\phi_-(t,1)
\end{bmatrix} = K\begin{bmatrix}
\phi_+(t,1)\\ \phi_-(t,0)
\end{bmatrix} +Bu(t) + w_b\\
&&\hspace*{-0.35cm}=
\begin{bmatrix}K_{11}&K_{12}\\K_{21}&K_{22}\end{bmatrix} \begin{bmatrix}
\phi_+(t,1)\\ \phi_-(t,0)
\end{bmatrix} + \begin{bmatrix}B_1\\B_2\end{bmatrix} u(t) + w_b\label{Hyp_BC_decomposed}
\end{eqnarray}
where  $\phi = \begin{bmatrix}
\phi_+\\ 
\phi_-
\end{bmatrix}$ with $\phi^+$ in $\RR^\ell$, $\phi^-$ in $\RR^{n-\ell}$ and where $w_b$ in $\RR^p$ is an unknown disturbance, 
 $u(t)$ is a control input taking values in $\RR^m$ and $K$, $B$  are matrices of appropriate dimensions.



The output to be regulated to a prescribed value denoted by $y_{ref}$, is given as a disturbed linear combination of the  boundary conditions. Namely, the outputs to regulate are in $\RR^m$ given as
\begin{equation}\label{eq_Output_old}
y(t) = L_1\begin{bmatrix}
\phi_+(t,0)\\
\phi_-(t,1)
\end{bmatrix} + L_2\begin{bmatrix}
\phi_+(t,1)\\ \phi_-(t,0)
\end{bmatrix} + w_y,
\end{equation}
where $L_1$ and $L_2$ are two matrices in $\RR^{m\times n}$ and $w_y$ is an unknown disturbance in $\RR^m$.
We wish to find a positive real number $k_i$ and a full rank matrix $K_i$ such that
\begin{equation}\label{eq_IntContr_Hyp}
u(t) = k_i K_iz(t)\ ,\ z_t(t)=y(t)-y_{ref} \ , z(0)=z_0
\end{equation}
where $z(t)$ takes value in $\RR^m$ and $z_0\in \RR^m$
solves the regulation problem  $\forall y_{ref}\in \RR^m$.

The state space denoted by $\XR_e$ 
of the system \eqref{eq_Hyp}-\eqref{Hyp_BC} in closed loop with the control law \eqref{eq_IntContr_Hyp}
is the Hilbert space defined as:
$$
\XR_e = (L^2(0,1),\RR^n)\times\RR^m,
$$
equipped with the norm defined for $ \varphi_e=(\phi,z)$ in $\XR_e$ as:
$$
\|v\|_{\XR_e} = \|\phi\|_{L^2((0,1),\RR^m)}  + |z|  .
$$
We introduce also a smoother state space defined as:
$$
\XR_{e1} = (H^1(0,1),\RR^n)\times\RR^m.
$$

\subsection{Output regulation result}
In this section, we give a set of sufficient conditions allowing  to solve the regulation problem as described in the introduction. 
Our approach follows what we have done in the former section.
Following \cite[Proposition 5.1, p161]{BastinCoron_2016stability} we consider the following assumption.

\begin{assumption}[Input-to-State Exponential Stability]\label{ass_HypEXPStab}
There exist a $C^1$ function   $P:[0,1]\rightarrow\mathcal D_{n}$, a real numbers $\mu>0$, $\underline P$, $\overline P$ and a positive definite matrix $S$ in $\RR^{n\times n}$ such that 
\begin{eqnarray}\label{eq_der_LyapHyp}
 (P(s) \Lambda_0(s) )_s 
-P(s)\Lambda_1(s) -\Lambda_1^\top(s)P(s)\nonumber \\ \leqslant   -\mu P(s),\\
 \label{eq_P}
\underline P \Id_n \leqslant   P(s) \leqslant   \overline P \Id_n\ ,\ \forall s\in[0,1],
\end{eqnarray}
and
\begin{equation}\label{eq_Lyap_Boundary}
-K_+^\top P(1)\Lambda_0(1)
K_+
+
K_-^\top P(0)\Lambda_0(0)K_- 
\leqslant   -S.
\end{equation}
where
\begin{equation}\label{eq_K}
K_+ = \begin{bmatrix}
\Id_\ell & 0\\
K_{21} & K_{22}
\end{bmatrix}\ ,\ K_- = \begin{bmatrix}
K_{11} & K_{12}\\
0 & \Id_{n-\ell}
\end{bmatrix}
\end{equation}
\end{assumption}
\vspace*{0.2cm}

As it will be seen in the following section, this assumption is a sufficient condition for exponential stability of the equilibrium of the open loop system.
It can be found in \cite{BastinCoron_2016stability} in the case in which  $S$ may be semi-definite positive. 
The positive definiteness of $S$ is fundamental to get an  input-to-state stability (ISS) property of the open loop system with respect to the disturbances on the boundary.
More general results are given in \cite{prieur2012iss}.

\noindent The second assumption is related to the rank condition.
Let $\Phi:[0,1]\rightarrow\RR^{n\times n}$ be the  matrix function solution to the system
$$
\begin{aligned}
\Phi_s(s) &= \Lambda_0(s)^{-1}\Lambda_1(s) \Phi(s),
\\
\Phi(0)&=\Id_n .
\end{aligned}
$$
We denote $\Phi(s) = \begin{bmatrix}
\Phi_{11}(s)&\Phi_{12}(s)\\
\Phi_{21}(s)&\Phi_{22}(s)
\end{bmatrix}$
and
$$
\Phi_+(1) = \begin{bmatrix}
\Phi_{11}(1) & \Phi_{12}(1)\\
0&\Id_{n-\ell}
\end{bmatrix}\ , \ 
\Phi_-(1) = \begin{bmatrix}
\Id_\ell  & 0\\
\Phi_{21}(1)&\Phi_{22}(1)
\end{bmatrix}
$$

\begin{assumption}[Rank condition 1]\label{ass_RankCondition}
The matrix in $\RR^{n\times n}$  $
\Phi_-(1)-K\Phi_+(1)
$ is full rank and so is the matrix $T$ defined as
\begin{equation}\label{eq_Ki_Hyp}
T_1 = \left(L_1\Phi_-(1)+L_2\Phi_+(1)\right)
\left(\Phi_-(1)-K\Phi_+(1)\right)^{-1}B.
\end{equation}
\end{assumption}

\vspace{0.5em}

Another rank condition has to be introduced. This one is used when solving the forwarding equation.
Let $\Psi:[0,1]\mapsto \RR^{n\times n}$ be the matrix function solution to the system
\begin{equation}\label{eq_Psi}
\begin{aligned}
\Psi_s(s) &= \Psi(s) \left( \Lambda_1(s)-\Lambda_{0s}(s)\right)\Lambda_0(s)^{-1},\\
\Psi(0)&=\Id_n.
\end{aligned}
\end{equation}
\begin{assumption}[Rank condition 2]\label{ass_MCondition}
The matrix in $\RR^{n\times n}$  
\begin{equation}\label{eq_Matrice}
\Psi(1)\Lambda_0(1)K_+-\Lambda_0(0)K_-
\end{equation}
is full rank and so is the matrix
$$
T_2 = -L_1B
+M\left(\Lambda_0(0)\begin{bmatrix}B_1\\0\end{bmatrix}
-\Psi(1)\Lambda_0(1)\begin{bmatrix}0\\B_2\end{bmatrix}\right)
$$
where
\begin{equation}\label{eq_MMatrix}
M = 
 (L_1K +L_2)\left(\Lambda_0(0)K_--\Psi(1)\Lambda_0(1)K_+\right)^{-1}.
\end{equation}
\end{assumption}

\noindent With these assumptions, the following result may be obtained.

\begin{theorem}[Regulation for hyperbolic PDE systems]\label{Theo_HypRegulation}
Assume that Assumptions \ref{ass_HypEXPStab}, \ref{ass_RankCondition} and \ref{ass_MCondition} are satisfied then with
$K_i = T_2^{-1}$
there exists $k_i^*>0$ such that for all $0<k_i<k_i^*$ the output regulation is obtained.
More precisely, for all $(w_b,w_y,y_{ref})$ in $\RR^p\times\RR^m\times \RR^m$, the following holds. 
\begin{enumerate}
\item For all $(\phi_0,z_0)$ in $\XR_e$ (resp. $\XR_{1e}$) 
which satisfies the boundary conditions \eqref{Hyp_BC} (resp. the $C^1$ compatibility condition), there exists a unique weak solution to
\eqref{eq_Hyp}-\eqref{Hyp_BC}-\eqref{eq_IntContr_Hyp}
 that we denote $v$ and
which belongs to $C^0([0,+\infty);\XR_e)$ (Respectively, strong solution in:
\begin{equation}\label{eq_SmoothProperties}
C^0 ([0,+\infty);\XR_{e1})\cap C^1([0,+\infty);\XR_e)\ ).
\end{equation}
\item There exists an equilibrium state denoted $v_\infty$ in $\XR_e$ which is globally exponentially stable in $\XR_e$ for system \eqref{eq_Hyp}-\eqref{Hyp_BC}-\eqref{eq_IntContr_Hyp}. More precisely, we have for all $t\geqslant   0$:
\begin{equation}\label{eq_ExpStabHyp}
\|v(t)-v_\infty\|_{\XR_e} \leqslant   k \exp(-\nu t) \|v_0-v_\infty\|_{\XR_e}.
\end{equation}
\item Moreover, if $v_0$ satisfies the $C^1$-compatibility condition and  is in $\XR_{1e}$, the regulation is achieved, i.e.
\begin{equation}\label{eq_RegulNorm}
\lim_{t\rightarrow +\infty} |y(t) - y_{ref}| = 0. 
\end{equation}
\end{enumerate}
\end{theorem}

\par\vspace{0.5em}
The next section is devoted to the proof of this result.

\subsection{About this result}

The first assumption needed in Theorem \ref{Theo_HypRegulation} is Assumption \ref{ass_HypEXPStab}. When considering only integral control laws, there is no hope to obtain the result without assuming exponential stability of the open loop system.
Assumption \ref{ass_HypEXPStab} is slightly more restrictive than exponential stability since it requires an ISS property with respect to the input $u$. 
In the case in which this assumption is not satisfied for a given hyperbolic system, a possibility is to modify the boundary condition via a static output feedback (or proportional feedback) following the route of \cite{BastinCoron_2016stability} in order  to satisfy this assumptions.

One interest of our approach is that, part of the exponential stability of the closed loop system, only
Assumptions \ref{ass_RankCondition} and \ref{ass_MCondition} which are rank conditions involving the boundary conditions have to be satisfied.
In the case in which the two above mentioned assumptions  are not satisfied, we may obtain these properties by adding a proportional feedback and consequently changing the value of $K$ in $T_1$ and $T_2$ to obtain these rank conditions.
These Assumptions \ref{ass_RankCondition} and \ref{ass_MCondition} are version of Point 3) in Theorem \ref{Theo_XuJerbi}. 

In the particular case in which $\Lambda_0$ is constant and $\Lambda_1=0$, the matrix function $\Phi(s)$ and $\Psi(s)$ are simply equal to identity for all $s$ in $[0,1]$. 
In that case, it yields
\begin{equation}\label{eq_T1Simple}
T_1 = (L_1 + L_2)(\Id_n-K)^{-1}B,
\end{equation}
and,
\begin{align}
T_2 &= -L_1 B + (L_1K + L_2)(K_- - K_+)^{-1}\begin{bmatrix}B_1\\-B_2\end{bmatrix}\nonumber\\
&=\begin{multlined}[t]
-L_1B + (L_1K+L_2)\begin{bmatrix}K_{11} - \Id_\ell & K_{12}\\
-K_{21}&\Id_{n-\ell}-K_{22}\end{bmatrix}^{-1}\\\times\begin{bmatrix}\Id_\ell&0\\0&-\Id_\ell\end{bmatrix}B
\end{multlined}\nonumber\\
&=\begin{multlined}[t]
-L_1B + (L_1K+L_2)\\
\times\left(\begin{bmatrix}\Id_\ell&0\\0&-\Id_\ell\end{bmatrix}\begin{bmatrix}K_{11} - \Id_\ell & K_{12}\\
-K_{21}&\Id_{n-\ell}-K_{22}\end{bmatrix}\right)^{-1}
B
\end{multlined}\nonumber\\
&=\begin{multlined}[t]
-L_1B - (L_1K+L_2)(\Id_n-K)^{-1}
B
\end{multlined}\nonumber\\
&=\begin{multlined}[t]
\left[-L_1(\Id_n-K) - (L_1K+L_2)\right](\Id_n-K)^{-1}
B
\end{multlined}\nonumber\label{eq_T2Simple}\\
&=-(L_1+L_2)(\Id_n-K)^{-1}B.
\end{align}
Hence when $\Lambda_0$ is constant and $\Lambda_1=0$, Assumption \ref{ass_RankCondition} and Assumption \ref{ass_MCondition} are equivalent.

Also, an interesting aspect of this Lyapunov approach is that explicit values of the supremum value of the gain $k_i^*$ may be given.
For instance, as in \cite{TrinhAndrieuXu_TAC_17} consider the very particular case of a transport equation.
In this case the system is simply
\begin{align*}
&\phi_t(s,t) + \phi_s(s,t) = 0, s\in(0,1), t\in[0,+\infty)\\
&\phi(t,0) = u(t) + w_b\\
&y(t) = \phi(1,t) + w_y
\end{align*}
We can apply Theorem \ref{Theo_HypRegulation} with $n=1$, $\Lambda_0(s)=-1$, $\Lambda_1(s)=0$, $K=0$, $B=1$, $L_1=0$, $L_2=1$.
This yields $\Psi(s)=1$, $\Psi(s)=1$, $T_1=1$, $T_2=-1$.
Hence, Assumptions  \ref{ass_RankCondition} and \ref{ass_MCondition} are satisfied.
Assumption \ref{ass_HypEXPStab} is satisfied for all $\mu>0$ with $P(s)=e^{-\mu s}$, $S=1$, $\overline P=1$, $\underline P=e^{-\mu}$.
In that case, employing theorem \ref{Theo_HypRegulation}, it yields that there exists $k_i^*>0$ such that for all $0<k_i<k_i^*$ 
with $u(t)=-k_i z$, $\dot z = y$, the output regulation is obtained and so the output converges asymptotically to zero.
Following the proof of Theorem \ref{Theo_HypRegulation}, equation (\ref{eq_kstar}) gives
$$
k^*_i=\sqrt{\mu e^{-\mu}}.
$$
This bound is better then the one obtained in \cite{TrinhAndrieuXu_TAC_17} for the linear transport equation (its maximal value is obtained for $\mu=1$ and is $\frac{1}{\sqrt{e}}$. 
Note however that similar to the bound of \cite{TrinhAndrieuXu_TAC_17}, the result obtained with our novel Lyapunov functional is far from the value we get following a frequency approach ($\frac{\pi}{2}$ in this case).
Recently in \cite{coron2018pi}, the Lyapunov functional obtained in \cite{TrinhAndrieuXu_TAC_17}  has been modified to reach this optimal value of the integral gain. A natural question for future research topic is to know if it is possible to modify the Lyapunov functional obtained in Theorem \ref{Theo_HypRegulation} following the methods of \cite{coron2018pi} to remove the conservatism.

\subsection{Illustration in a $2 \times 2$ hyperbolic system}

Theorem \ref{Theo_HypRegulation} generalizes many available results on output regulation via integral action for hyperbolic PDEs available in the literature.
For instance, the case of $2\times 2$ linear hyperbolic systems has been considered in \cite{TrinhAndrieuXu_IFACWC_17},  \cite{dos2008boundary},  (see also \cite[Section 2.2.4]{BastinCoron_2016stability}).
The case of cascade of such systems is also considered in \cite{TrinhAndrieuXu_Aut_18}.
Note also that in \cite{TerrandJeanneEtAl_TAC_2018}, this procedure is applied on a Drilling model which is composed of a hyperbolic PDE coupled with a linear ordinary differential equation.

In order to compare the way we improve existing results, the same example as in \cite{dos2008boundary} is considered. 
In this context, the linearized de Saint-Venant equations can be written in the form of \eqref{eq_Hyp}-\eqref{Hyp_BC}.
After normalization, one gets~: 
\begin{eqnarray}
\Lambda_0(s)=\begin{bmatrix}
c &0\\
0& -d
\end{bmatrix} \text{ and }\Lambda_1(s)=0_{2 \times 2}\ , \forall s 
\end{eqnarray}
where $c>0$ and $d>0$ and 
\begin{eqnarray}
K=\begin{bmatrix}
0 &k_0\\
k_1& 0
\end{bmatrix} \text{ and }B=\begin{bmatrix}
b_0&0\\
 0&b_1
\end{bmatrix},
\end{eqnarray}
with $b_0\neq 0$ and $b_1\neq 0$. For the system \eqref{eq_Hyp}-\eqref{Hyp_BC} with these parameters, it is shown in \cite{dos2008boundary} that the  output of dimension $m=2$ defined in \eqref{eq_Output_old} with
\begin{eqnarray}\label{eq_SaintVenantOutput}
L_1=\begin{bmatrix}
\frac{c}{c+d} &0\\
0& \frac{-1}{c+d} 
\end{bmatrix} \text{ and }L_2=\begin{bmatrix}
0& \frac{d}{c+d}\\
 \frac{1}{c+d}&0 
\end{bmatrix}
\end{eqnarray}
 can be regulated with an integral control law provided
 \begin{equation}\label{eq_CondDos}
     |k_0k_1|<1\ , \ |k_0|<1\ , \ |k_1|<\frac{c}{d}.
 \end{equation}
 
On another hand, employing  (\cite{coron2018pi}-\cite{dos2008boundary}), Assumptions \ref{ass_HypEXPStab} is satisfied assuming that $|k_0k_1|<1$.
Moreover, 
with equations \eqref{eq_T1Simple} and \eqref{eq_T2Simple}, it yields,
$$
T_1 = -T_2 = \frac{1}{c+d}\begin{bmatrix}
c & d \\1 &-1
\end{bmatrix}\begin{bmatrix}1&-k_0\\-k_1 & 1\end{bmatrix}^{-1}\begin{bmatrix}
b_0&0\\
 0&b_1
\end{bmatrix}.
$$
This matrix is well defined and  full rank if $|k_0k_1|<1$ and consequently Assumptions \ref{ass_RankCondition} and \ref{ass_MCondition} are always satisfied.
Hence, employing Theorem \ref{Theo_HypRegulation}, both outputs defined in \eqref{eq_SaintVenantOutput} can be regulated with an integral control law with the only assumption that $|k_0k_1|<1$.\\
Then 
\begin{eqnarray}
&&\hspace*{-1cm}K_i=T_2^{-1}\nonumber\\
&&=\vartheta
\begin{bmatrix}
b_1(1-k_0) & b_1 (d+ck_0)\\
 b_0(1-k_1) & b_0(c+dk_1)
\end{bmatrix}\\
&&\hspace*{-1cm}\text{with }\vartheta=\frac{-(c+d)^2(1-k_0k_1)^2b_0^{-1}b_1^{-1}}
{ [(1-k_0)(c+dk_1)+(1-k_1)(d+ck_0)]}
\end{eqnarray} and 
\begin{eqnarray}
 k_i^* = \frac{\sqrt{\mu \underline P}}{|M|\overline \Psi \sqrt{c \left|T_2^{-1}\right|}}\nonumber\\
\end{eqnarray}
$\mu$ is given in \cite{dos2008boundary}, and $\underline P$ the lower bound of the Lyapunov can be deduced easily from the expression of the Lyapunov function involved.\\
$M$ has been defined above, with $T_2$. As $ \Psi $ is the identity matrix, $\overline \Psi $ is 1.
Remark that in \cite{dos2008boundary}, $K_i$ is diagonal and here is full matrix.
Note that some other choices of $K_i$ are possible as long as
$$
T_2K_i + K_i^\top T_2^\top >0.
$$
 To conclude, a work is needed to transpose this approach to the global de Saint-Venant equations this is the aim of another paper.

\section{Proof of Theorem \ref{Theo_HypRegulation}}
\label{Sec_ProofTheoHyper}
The proof of Theorem \ref{Theo_HypRegulation} is divided into three steps.
In a first part, it is shown that with Assumption \ref{ass_RankCondition}, it can be shown that the closed loop system \eqref{eq_Hyp}-\eqref{Hyp_BC}-\eqref{eq_IntContr_Hyp} admits a steady state.
In a second step, it is established that the desired regulation is obtained provided the steady state is exponentially stable.
Finally, the construction of an appropriate Lyapunov functional is performed to show the exponential stability of the equilibrium.
 
\subsection{Stabilization implies regulation}

In this first subsection, we explicitly give the equilibrium state of the system \eqref{eq_Hyp}-\eqref{Hyp_BC}-\eqref{eq_IntContr_Hyp}.
We show also that if we assume that $k_i$ and $K_i$ are selected such that this equilibrium point is exponentially stable along the closed loop, then the regulation is achieved.

\subsubsection{Definition of the equilibrium}

The first step of the study is to exhibit equilibrium denoted $\phi_\infty,z_\infty$ of the disturbed hyperbolic PDE in closed loop with the boundary integral control (i.e. system \eqref{eq_Hyp}-\eqref{eq_IntContr_Hyp}).

We have the following proposition.
\begin{proposition}\label{Prop_WellPosednessEquilibria}
Assumption \ref{ass_RankCondition} is a necessary and sufficient condition for the existence of an equilibrium of the system \eqref{eq_Hyp}-\eqref{Hyp_BC}-\eqref{eq_IntContr_Hyp}.
Moreover, if Assumption  \ref{ass_RankCondition} holds then point 1) of Theorem \ref{Theo_HypRegulation} holds.
\end{proposition}
\if\cdc1

The proof of this result has been removed due to space limitation and will be available in the journal version of this paper.

We can introduce $\tilde \phi(x,t) = \phi(x,t) - \phi_\infty$ and $\tilde z(t) = z(t)-z_\infty$.
It can be checked that $\tilde \phi, \tilde z$ satisfies the following system:
\begin{equation}\label{eq_Hyp_Unperturbed}
\begin{array}{rl}
\tilde \phi_t &= \Lambda \tilde \phi_x \ ,\ \varphi\in(0,1),\\
z_t&=L_1\begin{bmatrix}
\tilde \phi_+(t,0)\\
\tilde \phi_-(t,1)
\end{bmatrix} + L_2\begin{bmatrix}
\tilde \phi_+(t,1)\\ \tilde \phi_-(t,0)
\end{bmatrix}
\end{array}\ ,\ t\in[0,+\infty),
\end{equation}
with the boundary conditions
\begin{align}\label{Hyp_BC_Undisturbed}
\begin{bmatrix}
\tilde \phi_+(t,0)\\
\tilde \phi_-(t,1)
\end{bmatrix} &= K\begin{bmatrix}
\tilde \phi_+(t,1)\\ \tilde \phi_-(t,0)
\end{bmatrix} +Bu(t),\\
u(t)&=k_iK_i\tilde z(t) .\label{eq_BoundCont}
\end{align}
\else
\begin{proof}
First of all, equilibria are such that
$$
\phi_{\infty s}(s) = -\Lambda_0(s)^{-1}\Lambda_1(s)\phi_\infty(s),
$$
for all $s$ in $[0,1]$.
Hence,
\begin{equation}\label{eq_Solphiinfty}
\phi_{\infty}(s) = \Phi(s)\phi_{\infty}(0).
\end{equation}
Hence,
\begin{align*}
\begin{bmatrix}
\phi_{\infty+}(0)\\
\phi_{\infty-}(1)
\end{bmatrix} &= \Phi_-(1)
\phi_{\infty }(0),\\
\begin{bmatrix}
\phi_{\infty+}(1)\\
\phi_{\infty-}(0)
\end{bmatrix} &= \Phi_+(1)\phi_{\infty }(0)
\end{align*}
Moreover, with $z_t=0$, we have
$$
  L_1\begin{bmatrix}
\phi_{\infty+}(0)\\
\phi_{\infty-}(1)
\end{bmatrix} + L_2\begin{bmatrix}
\phi_{\infty+}(1)\\ \phi_{\infty-}(0)
\end{bmatrix} =y_{ref}- w_y,
$$
Hence,
\begin{equation}\label{eq_y_equil}
  \left(L_1 \Phi_-(1)+ L_2\Phi_+(1)\right)
\phi_{\infty}(0)
 =y_{ref}- w_y,
\end{equation}

On another side, boundary conditions \eqref{Hyp_BC} gives
\begin{eqnarray}\label{Hyp_BC2}
\left(\Phi_-(1) - K\Phi_+(1)\right)\phi_\infty(0) =
B k_i K_i z_\infty + w_b
\end{eqnarray}
For all $w_y$ and $y_{ref}$ both in $\RR^m$, $w_b$ in $\RR^p$, by Assumption \ref{ass_RankCondition}
 and  since the matrix $K_i$ is full rank the former equation and \eqref{eq_y_equil} admit a unique solution $(z_\infty,\phi_{\infty}(0))$ given as
\begin{multline}
z_\infty = \frac{ K_i^{-1}}{k_i}T_1^{-1}(y_{ref}-w_y) \\
-\frac{ K_i^{-1}}{k_i}T_1^{-1} \left(L_1\Phi_-(1)+L_2\Phi_+(1)\right)\\
\times
\left(\Phi_-(1)-K\Phi_+(1)\right)^{-1} w_b 
\end{multline}
and,
\begin{multline}
\phi_\infty(0) =
 \left(\Phi_-(1)-K\Phi_+(1)\right)^{-1}k_iBK_iz_\infty \\+ \left(\Phi_-(1)-K\Phi_+(1)\right)^{-1}w_b.
\end{multline}
Finally, in that case, we can introduce $\tilde \phi(s,t) = \phi(s,t) - \phi_\infty(s)$ and $\tilde z(t) = z(t)-z_\infty$.
It can be checked that $\tilde \phi, \tilde z$ satisfies the following system:
\begin{align}\label{eq_Hyp_Unperturbed}
&\begin{multlined}[t]
\tilde \phi_t(s,t)+ \Lambda_0(s) \tilde \phi_s(s,t) +\Lambda_1(s)\tilde\phi(s,t)\\=0\ ,\ s\in(0,1),
\end{multlined}
\end{align}
\begin{equation}
\tilde z_t=L_1\begin{bmatrix}
\tilde \phi_+(t,0)\\
\tilde \phi_-(t,1)
\end{bmatrix} + L_2\begin{bmatrix}
\tilde \phi_+(t,1)\\ \tilde \phi_-(t,0)
\end{bmatrix}
\end{equation}
with the boundary conditions
\begin{align}\label{Hyp_BC_Undisturbed}
\begin{bmatrix}
\tilde \phi_+(t,0)\\
\tilde \phi_-(t,1)
\end{bmatrix} &= K\begin{bmatrix}
\tilde \phi_+(t,1)\\ \tilde \phi_-(t,0)
\end{bmatrix} +Bu(t),\\
u(t)&=k_iK_i\tilde z(t) .\label{eq_BoundCont}
\end{align}
As it is shown in \cite{BastinCoron_2016stability}, for each initial condition $\tilde v_0=(\tilde \phi_0,\tilde z_0)$ in $\XR_e$ which satisfies the boundary conditions \eqref{Hyp_BC}, there exists a unique weak solution that we denoted $\tilde v$ and
which belongs to $C^0([0,+\infty);\XR_e)$. Moreover, if the initial condition $\tilde v_0$ satisfies also the $C^1$-compatibility condition (see \cite{BastinCoron_2016stability} for more details) and lies in $\XR_{e1}$ then the solution lies in the set defined in \eqref{eq_SmoothProperties}.
\end{proof}
\fi
\subsubsection{Sufficient conditions for Regulation}
In the following, we  show that the regulation problem can be rephrased as a stabilization of the equilibrium state introduced previously.

\begin{proposition}\label{Prop_StabImplyRegul}
Assume Assumption \ref{ass_RankCondition} holds and that there exist a functional $V_e:\XR_e\rightarrow\RR_+$, and positive real numbers $\mu_e$ and $L_e$ such that:
\begin{equation}\label{eq_LyapEquationW2}
\frac{\|v_\infty-v\|_{\XR_e}^2}{L_e} \leqslant    V_e(v) \leqslant    L_e\|v_\infty-v\|_{\XR_e}^2.
\end{equation}
Assume moreover that for all $v_0$ in $\XR_e$ and all $t_0$ in $\RR_+$ such that the solution $v$ of system \eqref{eq_Hyp}-\eqref{Hyp_BC}-\eqref{eq_IntContr_Hyp} initiated from $v_0$ is $C^1$ at $t=t_0$, we have:
\begin{equation}\label{eq_LyapEquationW}
\dot V_e(t) \leqslant    -\mu_e V_e(t),
\end{equation}
where with a slight abuse of notation $V_e(t) = V_e(v(t))$.
Then points 1), 2) and 3) of Theorem \ref{Theo_HypRegulation} hold.
\end{proposition}
\if\cdc1

\noindent\textit{Sketch of the proof:}
Point 1) is directly obtained from Proposition \ref{Prop_WellPosednessEquilibria}.
The proof of point 2) is by now standard.
The last point is obtained employing Sobolev embedding and showing the convergence of trajectories toward the equilibrium in the supremum norm. 
\vspace{0.5em}

\else
\begin{proof}
Point 1) is directly obtained from Proposition \ref{Prop_WellPosednessEquilibria}.
The proof of point 2) is by now standard. 
Let $v_0$ be in $\XR_{e1}$ and satisfies the $C^0$ and $C^1$-compatibility conditions. It yields that $v$ is $C^1$ for all $t$.
Consequently, \eqref{eq_LyapEquationW} is satisfied for all $t\geqslant   0$.
With Gr\"onwall's lemma, this implies that:\\
$
\hspace*{2cm}  V_e(v(t)) \leqslant    e^{-\mu_e t}V_e(v_0).\\
$
Hence with \eqref{eq_LyapEquationW2}, this implies that \eqref{eq_ExpStabHyp} holds with $k=L_e$ and $\nu=\frac{\mu_e}{2}$  for initial conditions in $\XR_{e1}$. $\XR_{e1}$ being dense in $\XR_e$, the result holds also with initial condition in $\XR_e$ and point 2) is satisfied.

On another hand, we have
\begin{align}\label{eq_outputFron}
y(t)-y_{ref} &=  L_1\begin{bmatrix}
\phi_+(t,0)\\
\phi_-(t,1)
\end{bmatrix} + L_2\begin{bmatrix}
\phi_+(t,1)\\ \phi_-(t,0)
\end{bmatrix} + w_y - y_{ref},\\
&=  L_1\begin{bmatrix}
\tilde \phi_+(t,0)\\
\tilde \phi_-(t,1)
\end{bmatrix} + L_2\begin{bmatrix}
\tilde \phi_+(t,1)\\ \tilde \phi_-(t,0)
\end{bmatrix} ,
\end{align}
with $\tilde \phi(t,x) = \phi(t,x)-\phi_\infty$.
To show that equation \eqref{eq_RegulNorm} holds, we need to show that the right hand side of the former equation tends to zero.  This may be obtained provided the initial condition is in  $\XR_1$.
Indeed, let $v_0$ be  
in $\XR_1$ and satisfies $C^1$-compatibility conditions. With (\ref{eq_SmoothProperties}), we know that $v_t \in C([0,\infty);\XR_e)$.
Moreover, 
$v_t$ satisfies the dynamics system \eqref{eq_Hyp}-\eqref{Hyp_BC}-\eqref{eq_IntContr_Hyp} with $w_b=0$, $w_y=0$, $y_{ref=0}$ (simply differentiate with time these equations). 
Hence, $\|v_t(t)\|_{\XR_e}$ converges exponentially toward $0$ and in particular
$$
\|\tilde \phi_{t}(t,\cdot)\|_{(L^2(0,1),\RR^n)}  \leqslant   k e^{-\nu t}\|v_t(0)\|.
$$
On another hand, employing \eqref{eq_Hyp}, it yields:
\begin{multline*}
\|\tilde \phi_{s}(t,\cdot)\|_{L^2((0,1),\RR^n)}  \\=\|\Lambda_0^{-1}(\tilde \phi_{t}(t,\cdot)+\Lambda_1(\cdot)\phi(t,\cdot))\|_{L^2((0,1),\RR^n)}
 .    
\end{multline*}
Hence,
\begin{multline}
 \|\tilde \phi_{s}(t,\cdot)\|_{L^2((0,1),\RR^n)} \\\leqslant    c\left(\|\tilde \phi_{t}(t,\cdot)\|_{L^2((0,1),\RR^n)}+\|\tilde \phi(t,\cdot))\|_{L^2((0,1),\RR^n)}\right)
 .
\end{multline}
where $c$ is a positive constant.
Consequently $\|\tilde \phi_s(t,\cdot)\|_{L^2((0,1),\RR^n)}$ converges also to zero and so is $\|\tilde \phi(t,\cdot)\|_{H^1((0,1),\RR^n)}$.
With Sobolev embedding
$$
\sup_{x\in[0,1]} |\tilde \phi(t,x)|\leqslant    C\|\tilde \phi(t,\cdot)\|_{H^1((0,1),\RR^n)},
$$
where $C$ is a positive real number. It implies that:
$$
\lim_{t\rightarrow+\infty}|\tilde \phi(t,1)|+|\tilde \phi(t,0)|=0.
$$
Consequently, with \eqref{eq_outputFron}, it yields that \eqref{eq_RegulNorm} holds and point 3) is satisfied.
\end{proof}
\fi
\noindent With this proposition in hand,  to prove the Theorem \ref{Theo_HypRegulation}, it is sufficient to construct a Lyapunov functional $V_e$ which satisfies \eqref{eq_LyapEquationW2}-\eqref{eq_LyapEquationW} along $C^1$-solutions of \eqref{eq_Hyp}-\eqref{Hyp_BC}-\eqref{eq_IntContr_Hyp} or equivalently along $C^1$-solutions of \eqref{eq_Hyp_Unperturbed}-\eqref{Hyp_BC_Undisturbed}.
This is considered in the next section following the route of Section~\ref{Sec_LyapAbstract}.

\subsection{Lyapunov functional construction}

\subsubsection{Open loop ISS}
Inspired by the Lyapunov functional construction introduced in \cite{CoronBastinAndreaNovel_2008dissipative} (see also \cite{BastinCoron_2016stability}), we
know that typical Lyapunov functionals allowing to exhibit stability property for this type of hyperbolic PDE are given as
 functional $V:L^2((0,1),\RR^n)\rightarrow\RR_+$ defined as
\begin{equation}
V(\varphi) = \int_0^1 \varphi(s)^\top P(s) \varphi(s) ds\ ,
\end{equation}
where $P:[0,1]\rightarrow\mathcal D_{n}$ is a $C^1$ function. Typically in \cite{CoronBastinAndreaNovel_2008dissipative}, these functions are taken as exponential.

With a slight abuse of notation, we write $V(t)=V(\tilde \phi(\cdot,t))$ and we denote by $\dot V(t)$ the time derivative of the Lyapunov functional along solutions which are $C^1$ in time. 
In our context, with Assumption \ref{ass_HypEXPStab}, it yields the following proposition.
\begin{proposition}
If Assumption \ref{ass_HypEXPStab} holds, 
there exists a positive real number $c$  such that for every
solution $\phi$ of \eqref{eq_Hyp}-\eqref{Hyp_BC} initiated from $(\tilde \phi_0,\tilde z_0)$ in $\XR_e$ which satisfies \eqref{Hyp_BC_Undisturbed} 
\begin{equation}\label{eq_OpenLoopLyapHyp}
\dot V(t) \leqslant   -\mu V(t) + c|u(t)|^2\ .
\end{equation}
\end{proposition}
\if\cdc1
\vspace{0.5em}
The proof of this proposition has been removed due to space limitation.
\else
\begin{proof}
First of all, with \eqref{eq_Hyp}, 
\begin{multline*}
\dot V(t) 
= - \int_0^1 2\phi(t,s)^\top P(s) \Lambda_0(s)\phi_s(t,s)ds \\ -\int^1_0 \phi(t,s)^\top \left(P(s) \Lambda_1(s) + \Lambda_1(s)^\top P(s)\right)\phi(t,s) 
ds\\
\end{multline*}
With an integration by part, this implies
\begin{multline*}
\dot V(t) 
=\int_0^1 \phi(t,s)^\top \left[ (P(s) \Lambda_0(s))_s \right]\phi(t,s) ds\\
-\int_0^1 \phi(t,s)^\top (P(s) \Lambda_1(s) + \Lambda_1(s)^\top P(s))\phi(t,s) ds\\-
\phi(t,1)^\top P(1)\Lambda_0(1)\phi(t,1) \\+
\phi(t,0)^\top P(0)\Lambda_0(0)\phi(t,0).
\end{multline*}
With \eqref{eq_der_LyapHyp}, it gives
\begin{multline*}
\dot V(t) 
\leqslant  -\mu V(t)-
\phi(t,1)^\top P(1)\Lambda_0(1)\phi(t,1) \\+
\phi(t,0)^\top P(0)\Lambda_0(0)\phi(t,0).
\end{multline*}
With the boundary condition \eqref{Hyp_BC_decomposed} and \eqref{eq_Lyap_Boundary}, this implies
\begin{multline}\label{eq_dV}
\dot V(t) 
\leqslant  -\mu V(t) - \begin{bmatrix}
\phi_+(1)^\top &
\phi_-(0)^\top\end{bmatrix}S \begin{bmatrix}
\phi_+(1) \\
\phi_-(0)\end{bmatrix}\\
+ 2\begin{bmatrix}
\phi_+(1)^\top &
\phi_-(0)^\top\end{bmatrix}Q u(t) + u(t)^\top R u(t)
,
\end{multline}
where,
\begin{multline}
R = -\begin{bmatrix}
0 &
B_2^\top
\end{bmatrix}(P(1)\Lambda_0(1) + \Lambda_0(1) P(1))
\begin{bmatrix}
0 \\
B_2
\end{bmatrix} \\+ 
\begin{bmatrix}
B_1^\top &
0
\end{bmatrix}(P(0)\Lambda_0(0) + \Lambda_0(0) P(0))
\begin{bmatrix}
0 \\
B_1
\end{bmatrix} ,
\end{multline}
and,
\begin{multline*}
    Q = -K_+^\top(P(1)\Lambda_0(1) + \Lambda_0(1) P(1))\begin{bmatrix}
0 \\
B_2
\end{bmatrix}\\
+ K_-^\top (P(0)\Lambda_0(0) + \Lambda_0(0) P(0))\begin{bmatrix}
B_1 \\
0
\end{bmatrix}.
\end{multline*}
Since $S$ is positive definite, selecting $c$ sufficiently large, it yields 
$$
\begin{bmatrix}
-S & Q\\Q^\top & R-c\Id_m
\end{bmatrix}\leqslant   0\ .
$$
Consequently, \eqref{eq_dV} implies that \eqref{eq_OpenLoopLyapHyp} holds.
\end{proof}


\subsubsection{Forwarding approach to deal with the integral part}
Following the route of Section \ref{Sec_LyapAbstract}, a Lyapunov functional is designed from $V$ adding some terms to take into account the state of the integral controller.
Let the operator $\MR:L^1((0,1);\RR^n) \rightarrow \RR^m$ be given as
\begin{align}
\mathcal M\varphi &=\int_0^1M\Psi(s)\varphi(s)ds
\end{align}
where $\Psi$ is the matrix function defined in \eqref{eq_Psi},
and $M$
 is a matrix in $\RR^{m\times n}$ defined in \eqref{eq_MMatrix}.

Following the Lyapunov functional construction in Theorem \ref{Theo_Abstract}, we consider the candidate Lyapunov functional $V_e:L^2((0,1);\RR^n)\times\RR^m$ given as
\begin{equation}
V_e(\varphi,z) = V(\varphi) +  p (z-\mathcal \MR \varphi)^\top (z-\mathcal \MR \varphi).
\end{equation}
In the following theorem, it is shown that by selecting properly $K_i$, $k_i$ and $p$, this function is indeed a Lyapunov functional for the closed loop system.
Again, with a slight abuse of notation, we write $V_e(t)=V_e(\tilde \phi(\cdot,t),\tilde z(t))$ and we denote by $\dot V_e(t)$ the time derivative of the Lyapunov functional along solutions which are $C^1$ in time. 
\begin{proposition}
\label{Prop_Lint}
Assume that Assumptions \ref{ass_HypEXPStab} and \ref{ass_RankCondition} hold.
Then there exists a matrix $K_i$ in $\RR^{m\times m}$ and $k_i^*>0$ such that for all $0<k_i< k_i^*$,
there exist positive real numbers 
$L_e$ and $\mu_e$ such that for all $(\varphi,z)$ in $\XR_e$
\begin{equation}\label{eq_BoundW}
\frac{1}{L_e}\left(\|\varphi\|_\XR^2 +|z|^2\right)\leqslant   V_e(\varphi,z) \leqslant   L_e\left(\|\varphi\|_\XR^2 +|z|^2\right),
\end{equation}
and
along $C^1$ solution of the system \eqref{eq_Hyp_Unperturbed}-\eqref{Hyp_BC_Undisturbed}-\eqref{eq_BoundCont}
\begin{equation}\label{eq_W_NoPsi}
\dot V_e(t) \leqslant    -\mu_{e} V_e(t) \ ,\ \forall t\in\RR_+.
\end{equation}
\end{proposition}
\if\cdc1

The proof of this proposition is removed due to space limitation.

\else
\begin{proof}
With \eqref{eq_P}, it yields for all $\varphi$ in $L^2((0,1);\RR^n)$,
\begin{equation}\label{eq_BoundV}
\underline P  \|\varphi\|_{L^2((0,1);\RR^n)}^2 \leqslant   V(\varphi) \leqslant   \overline P \|\varphi\|_{L^2((0,1);\RR^n)}^2.
\end{equation}

Let $\overline \Psi>0$ be such that
$$
|\Psi(s)| \leqslant   \overline \Psi\ ,\ \forall s\in[0,1].
$$
Note that for all $\varphi$ in $L^2((0,1);\RR^n$, by Cauchy-Schwartz inequality,
\begin{equation}\label{eq_BoundM}
|\MR \varphi| \leqslant   |M|\, {\overline \Psi}\, \|\varphi\|_{L^2((0,1);\RR^n)} .
\end{equation}
Hence, for each $p>0$ equation \eqref{eq_BoundW} holds.

Note that along $C_1$ solution to system \eqref{eq_Hyp_Unperturbed}, we have
\begin{align*}
\mathcal M \tilde \phi_t(t,\cdot) &= \mathcal M(-\Lambda_0(\cdot) \tilde \phi_s(t,\cdot)-\Lambda_1(\cdot)\tilde \phi(t,\cdot)) \\
&= \int_0^1M\Psi(s)(-\Lambda_0(s) \tilde \phi_s(t,s)-\Lambda_1(s)\tilde \phi(t,s))ds
\end{align*}
With an integration by part this implies
\begin{multline*}
\mathcal M \tilde \phi_t(t,\cdot) 
= \int_0^1M(\Psi(s)\Lambda_0(s))_s \tilde \phi(t,s)ds\\
-\int_0^1M\Psi(s)\Lambda_1(s)\tilde \phi(t,s))ds\\
-M\left(\Psi(1)\Lambda_0(1)\tilde \phi(t,1)-\Lambda_0(0)\tilde \phi(t,0) \right).
\end{multline*}
This gives,
\begin{multline*}
\mathcal M \tilde \phi_t(t,\cdot) 
=\\ \int_0^1M\left(\Psi_s(s)\Lambda_0(s) +  \Psi(s)(\Lambda_{0s}(s)-\Lambda_1(s))\right)\tilde \phi(t,s)ds
\\
-M\left(\Psi(1)\Lambda_0(1)\tilde \phi(t,1)-\Lambda_0(0)\tilde \phi(t,0) \right).
\end{multline*}
With the definition of $\Psi$, it yields
$$
\mathcal M \tilde \phi_t(t,\cdot) 
=
-M\left(\Psi(1)\Lambda_0(1)\tilde \phi(t,1)-\Lambda_0(0)\tilde \phi(t,0) \right).
$$

With the boundary condition \eqref{Hyp_BC_Undisturbed}, it yields
\begin{eqnarray}
&&\hspace*{-1cm}\mathcal M \tilde \phi_t(t,\cdot) 
= \nonumber
-M\left(\Psi(1)\Lambda_0(1)K_+-\Lambda_0(0)K_- \right)\begin{bmatrix}
\tilde \phi_+(t,1)\\\tilde\phi_-(t,0)
\end{bmatrix}\\
&&\hspace*{-0.1cm}-M\Psi(1)\Lambda_0(1)\begin{bmatrix}0\\B_2\end{bmatrix}u(t)\nonumber
+M\Lambda_0(0)\begin{bmatrix}B_1\\0\end{bmatrix}u(t)
\end{eqnarray}
Hence,
with the definition of $M$, it implies
\begin{multline*}
\mathcal M \tilde \phi_t(t,\cdot) 
=
(L_1K+L_2)\begin{bmatrix}
\tilde \phi_+(t,1)\\\tilde\phi_-(t,0)
\end{bmatrix}\\
+M\left(\Lambda_0(0)\begin{bmatrix}B_1\\0\end{bmatrix}
-\Psi(1)\Lambda_0(1)\begin{bmatrix}0\\B_2\end{bmatrix}\right)u(t)
\end{multline*}

On another hand,
$$
z_t(t) = (L_1K+L_2)\begin{bmatrix}
\tilde\phi_+(t,1)\\\tilde \phi_-(t,0)
\end{bmatrix}+L_1Bu(t).
$$
This gives,
\begin{multline*}
 \mathcal M \tilde \phi_t(t,\cdot) 
  =z_t(t)- L_1Bu(t)\\+M\left(\Lambda_0(0)\begin{bmatrix}B_1\\0\end{bmatrix}
-\Psi(1)\Lambda_0(1)\begin{bmatrix}0\\B_2\end{bmatrix}\right)u(t).
\end{multline*}
Hence, it yields
\begin{equation}
 \mathcal M \tilde \phi_t(t,\cdot) 
  =z_t(t)+T_2u(t).
\end{equation}
We recognize here equation \eqref{eq_M} when $u=0$.
This gives with \eqref{eq_OpenLoopLyapHyp}
\begin{multline}
\dot V_e(t) 
\leqslant  -\mu V(t)+ c|u(t)|^2  \\
- 2p(z(t)-\mathcal M\phi(\cdot,t))^\top T_2 u(t).
\end{multline}
Let now, $K_i = T_2^{-1}$.
Hence, this gives with $u=k_iK_iz$,
\begin{multline}
\dot V_e(t)    \leqslant   -\mu V(t)  + ck_i^2 |K_iz(t)|^2\\
- 2p |z(t)|^2 k_i+ 2pk_i\mathcal (\MR\phi(\cdot,t))^\top z(t),
\end{multline}
With \eqref{eq_BoundM}, and completing the square it yields for all $\varphi$ in $L^2((0,1);\RR^n)$ and $z$ in $\RR^m$,
\begin{align}
2(\MR\varphi)^\top z &\leqslant   |\MR\varphi|^2 + |z|^2,\\
&\leqslant   |M|^2\overline \Psi ^2\|\varphi \|_{L^2((0,1);\RR^n)}^2 + |z|^2.
\end{align}
Merging the last two inequality yields,
\begin{multline}
\dot V_e(t)   \leqslant   -\mu V(t)+ pk_i|M|^2\overline \Psi ^2\|\phi(\cdot,t) \|_{L^2((0,1);\RR^n)}^2   \\
+ \left( ck_i^2 |K_i|^2-pk_i\right)  |z(t)|^2.
\end{multline}
With \eqref{eq_BoundV}, this yields
\begin{multline}
\dot V_e(t)   \leqslant   \left(-\mu +pk_i\frac{|M|^2\overline \Psi ^2}{\underline P}\right) V(t) \\
+  \left( ck_i^2 |K_i|^2-pk_i\right)  |z(t)|^2.
\end{multline}
Note that if
$$
pk_i < \frac{\mu \underline P}{|M|^2\overline \Psi^2  }\ ,\ k_i^2 < \frac{pk_i}{c \left|T_2^{-1}\right|},
$$
this yields the existence of $\mu_{e}$ such that equation (\ref{eq_W_NoPsi}) holds.
This is obtained for all $k_i<k_i^*$ when
\begin{equation}\label{eq_kstar}
    k_i^* = \frac{\sqrt{\mu \underline P}}{|M|\overline \Psi \sqrt{c \left|T_2^{-1}\right|}},
\end{equation}
and 
$$
p < \frac{\mu \underline P}{k_i|M|^2\overline \Psi^2 }.
$$
\end{proof}
\fi

With this proposition, the proof of Theorem \ref{Theo_HypRegulation} is completed.

\section{conclusion}
In the last three decades, the regulation problem has been studied for different classes of distributed parameter systems. 
Most of existing results follow a semigroup approach and the perturbation theory for linear operator. 
In this paper we have shown that is was also possible to construct Lyapunov functionals to address the regulation problem in the case in which is used an integral action.
This framework allows to explicitly give an integral gain.
Moreover, it is no more necessary to impose boundedness of control or measurement operators to guarantee the regulation.
This is applied to PDE hyperbolic systems and this allows to generalize many available results in this field.

\bibliographystyle{plain}

\bibliography{BibVA}

\begin{IEEEbiography}[{\includegraphics[width=1in,height=1.0in,clip,keepaspectratio]{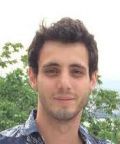}}]
{Alexandre Terrand-Jeanne} 
graduated in electrical engineering from ENS Cachan, France, in 2013. After one year in the robotic laboratory "Centro E.Piaggio" in Pisa, Italy, he is currently a doctoral student in LAGEP, university of Lyon 1. His PhD topic concerns the stability analysis and control laws design for systems involving hyperbolic partial differential equations coupled with nonlinear ordinary differential equations. This work is under the supervision of V. Dos Santos Martins, V. Andrieu and M. Tayakout-Fayolle. 
\end{IEEEbiography}

\begin{IEEEbiography}[{\includegraphics[width=1in,height=1.0in,clip,keepaspectratio]{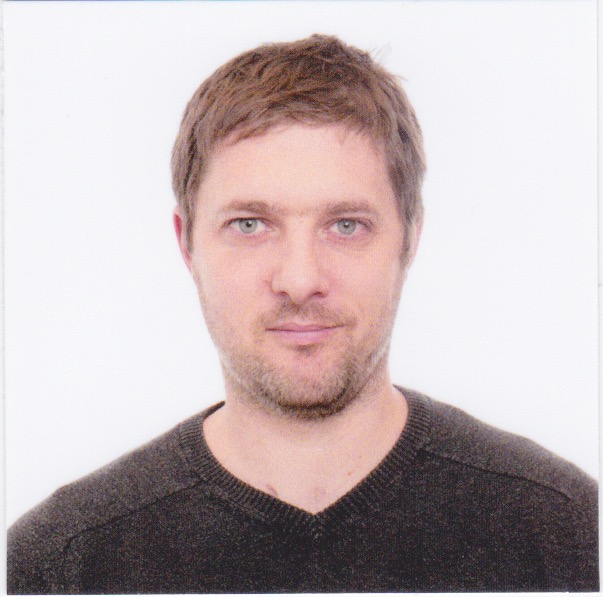}}]
{Vincent Andrieu}
graduated in applied mathematics from INSA de Rouen,
France, in 2001. After working in ONERA (French aerospace research company),
he obtained a PhD degree from Ecole des Mines de Paris in 2005. In 2006, he
had a research appointment at the Control and Power Group, Dept. EEE, Imperial
College London. In 2008, he joined the CNRS-LAAS lab in Toulouse, France, as
a CNRS-charg\'{e} de recherche. Since 2010, he has been working in LAGEP-CNRS,
Universit\'{e} de Lyon 1, France. In 2014, he joined the functional analysis
group from Bergische Universit\"at Wuppertal in Germany, for two sabbatical years.
His main research interests are in the feedback stabilization of controlled
dynamical nonlinear systems and state estimation problems. He is also interested
in practical application of these theoretical problems, and especially in the
field of aeronautics and chemical engineering.
Since 2018 he is an associate editor of the IEEE Transactions on Automatic Control, System \& Control Letters and IEEE Control Systems Letters.
\end{IEEEbiography}

\begin{IEEEbiography}[{\includegraphics[width=1in,height=1.0in,clip,keepaspectratio]{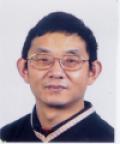}}]
{Cheng-Zhong XU}  received the Ph.D. degree in automatic control and signal processing from
Institut National Polytechnique de Grenoble, Grenoble, France, in 1989, and the Habilitation
degree in applied mathematics and automatic control from University of Metz, Metz, France,
in 1997. From 1991 to 2002, he was a Charg\'e de Recherche (Research Officer) in the Institut
National de Recherche en Informatique et en Automatique.
Since 2002, he has been a Professor of automatic control at the University of Lyon, Lyon, France.
His research interests include control of distributed parameter systems and its applications to
mechanical and chemical engineering.
He was an associated editor of the IEEE Transactions on Automatic Control, from 1995 to 1998.
He was an associated editor of the SIAM Journal on Control and Optimization, from 2011 to
2017.
\end{IEEEbiography}

\begin{IEEEbiography}[{\includegraphics[width=1in,height=1.0in,clip,keepaspectratio]{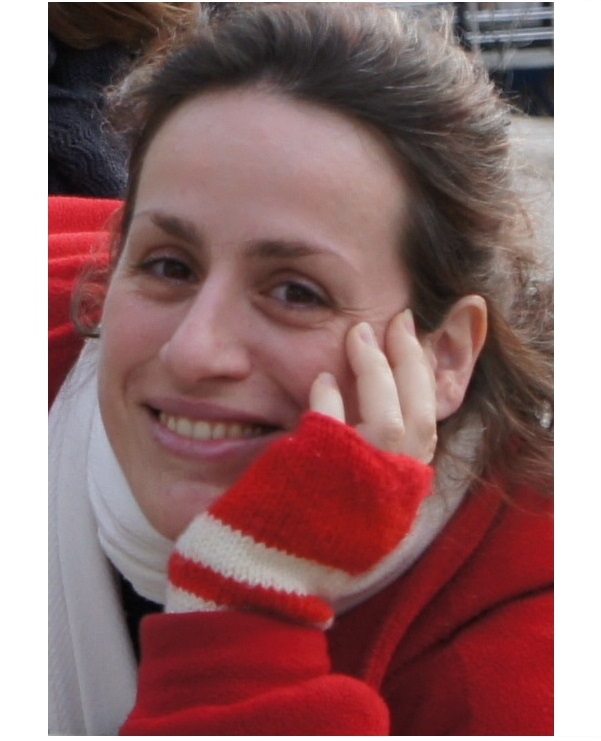}}]
{Val\'erie Dos Santos Martins} graduated in Mathematics from the University of Orl\'eans, France in 2001. She
received the Ph.D degree in 2004 in Applied Mathematics from the
University of Orl\'eans. After one year in the laboratory of
Mathematics MAPMO in Orl\'eans as ATER, she was post-doct in the
laboratory CESAME/INMA of the  University Catholic of Louvain, Belgium. Currently, she is professor assistant in the  laboratory LAGEP, University of Lyon 1. Her current research interests
include nonlinear control theory, perturbations theory of operators and semigroup, spectral theory
 and control of nonlinear partial differential equations.
\end{IEEEbiography}
\end{document}